\documentclass[11pt, twoside, a4paper]{article}
\usepackage{makeidx}\usepackage{amssymb}

\newcommand{\prs}{\langle\;,\;\rangle}
\newcommand{\too}{\longrightarrow}

\newcommand{\A}{{\cal A}}

\newcommand{\G}{{\mathfrak{g}}}

\newcommand{\h}{{\cal H}}

\newcommand{\B}{{\cal B}}

\newcommand{\D}{{\cal D}}

\newcommand{\di}{\displaystyle}

\newcommand{\al}{\alpha}
\newcommand{\be}{\beta}

\newcommand{\la}{\lambda}

\newcommand{\de}{\delta}

\font\bb=msbm10

\def\B{\hbox{\bb B}}
\def\R{\hbox{\bb R}}

\newtheorem{theo}{Theorem}[section]
\newtheorem{pr}{Proposition}[section]
\newtheorem{Le}{Lemma}[section]

\newtheorem{exem}{Example}
\newtheorem{rem}{Remark}

\title{ Ricci flat left invariant Lorentzian metrics on 2-step nilpotent Lie groups}

\author{Mohamed Boucetta}\date{}

\begin{document}\maketitle

\begin{abstract}
We determine all Ricci flat left invariant Lorentzian metrics on simply connected 2-step nilpotent Lie
groups. We show that the $2k+1$-dimensional Heisenberg Lie group $H_{2k+1}$ carries a Ricci
flat left invariant Lorentzian metric if and only if $k=1$. We show also that for any $2\leq q\leq k$, $H_{2k+1}$
carries a Ricci flat left invariant pseudo-Riemannian metric of signature $(q,2k+1-q)$ and we give explicite examples of such metrics.

\end{abstract}

{\it 2000 Mathematical Subject Classification: 53C50; Secondary
22E60, 53B30.

}

\section{Introduction}
Since Milnor's survey article \cite{milnor} which has already become a classic reference, the geometric properties of
Lie groups with left invariant Riemannian metrics have been studied extensively by many authors (see for instance
\cite{ lionel, dotti}). In \cite{milnor}, Milnor showed that   a Lie group carries a flat left invariant Riemannian
metric if and only if its Lie algebra  is a semi-direct product of an abelian algebra $\mathfrak{b}$ with an abelian
ideal $\mathfrak{u}$ and, for any $u\in\mathfrak{b}$, $ad_u$ is skew-symmetric. On the other hand,
a Ricci flat left invariant Riemannian metric must be flat (see \cite{besse}).
However, only a few partial results in the line of Milnor's study were known for pseudo-Riemannian left invariant
metrics. For instance,
 there exists some partial results on flat left invariant pseudo-Riemannian metrics on Lie groups
 (see \cite{guediri2, medina}) and there exists  Ricci flat left invariant pseudo-Riemannian metrics which are
 not flat (see \cite{jan}). Thus the study of flat or Ricci flat left invariant pseudo-Riemannian metrics is an
 open problem. In this paper, we study Ricci flat left invariant Lorentzian metrics on 2-step nilpotent Lie groups. We restrict our self to these groups because, first of all,  every thing in these groups is explicitly calculable and because of the richness of their geometry. Indeed, the geometry of left invariant Riemannian or Lorentzian metrics on 2-step nilpotent Lie groups were studied extensively by many authors (see for instance \cite{cordero, Eberlein, guediri1, guediri2}), partly because their relevance to General Relativity where they can be used to provide interesting counter-examples. Moreover,  Ricci flat Lorentzian manifolds are vacuum solutions of Einstein's equation and, as an interesting application of our main result, one  can built (by considering  quotients by lattices) a large class of compact Ricci flat Lorentzian manifolds. Note that all these Lorentzian manifolds are complete given that all left invariant pseudo-Riemannian metrics on a 2-step nilpotent Lie group are geodesically complete (see \cite{guediri3}).

 Let us give a brief outline of the results of this paper. In Section \ref{section1}, we give some properties of the Lie algebra of skew-symmetric endomorphisms on a pseudo-Euclidean vector space. In Section \ref{section2}, we establish a key formula giving the Ricci curvature of a pseudo-Euclidean 2-step nilpotent algebra (Lemma \ref{mainlemma}) and we give some of its immediate consequences. This formula is simple and plays a crucial role in this paper. In Section \ref{section3}, we show that the $2k+1$-dimensional Heisenberg Lie group $H_{2k+1}$ carries a Ricci flat left invariant Lorentzian metric if and only if $k=1$ and,  for any $2\leq q\leq k$, $H_{2k+1}$ carries a Ricci flat left invariant pseudo-Riemannian metric of signature $(q,2k+1-q)$ (Theorem \ref{heisenberg}), some explicite examples will be given (see Example \ref{example1}). In Section \ref{section4},
we determine all Ricci flat left invariant Lorentzian metrics on simply connected 2-step nilpotent Lie
groups (Theorems \ref{lorentzflat}-\ref{main}).
\section{Preliminaries}\label{section1}
In this section, we  give some properties of the Lie algebra of skew-symmetric endomorphisms on a pseudo-Euclidean
vector space.

A \emph{pseudo-Euclidean  vector space } is  a real vector space  of finite dimension $n$ endowed with  a
nondegenerate inner product  of signature $(q,n-q)=(-\ldots-,+\ldots+)$.  When the signature is $(0,n)$
(resp. $(1,n-1)$) the space is called \emph{Euclidean} (resp. \emph{Lorentzian}). Through this paper, we suppose
$q\leq n-q$.

Let $(V,\prs)$ be a pseudo-Euclidean vector space whose signature is $(q,n-q)$. A family $(u_1,\ldots,u_s)$ of vectors in $V$ is called \emph{orthogonal}  if, for $i,j=1,\ldots,s$ and $i\not=j$, $\langle u_i,u_j\rangle=0$. This family is called \emph{orthonormal} if it is orthogonal and, for any $i=1,\ldots,s$, $\langle u_i,u_i\rangle=1$.\\
A \emph{pseudo-Euclidean basis} of $V$ is a basis $(e_1,\bar e_2,\ldots,e_q,\bar e_q,f_1,\ldots,f_{n-2q})$ satisfying, for $i,j=1,\ldots,q$, $k=1,\ldots,n-2q,$
$$\langle e_i,e_j\rangle=\langle \bar e_i,\bar e_j\rangle=\langle  e_i, f_k\rangle=\langle \bar e_i,f_k\rangle=0,\quad
\langle e_i,\bar e_j\rangle=\de_{ij},$$and
$(f_1,\ldots,f_{n-2q})$ is  orthonormal. When $V$ is Lorentzian, we call such a basis \emph{Lorentzian}.   Pseudo-Euclidean basis always exist.\\
We denote by $Sym^-(V)$ the space of skew-symmetric endomorphisms of $V$ and we define a product $\prs^*$ on $Sym^-(V)$ by putting
$$\langle J,K\rangle^*=-{\mathrm tr}(J\circ K),$$where ${\mathrm tr}$ denotes the trace. This product is nondegenerate and defines a pseudo-Euclidean product on $Sym^-(V)$.

It is well-known that if $(V,\prs)$ is Euclidean then $\prs^*$ is definite positive and, for any $J\in Sym^-(V)$, there exists an orthonormal basis $\B$ of $V$ and a family of real numbers  $0<\la_1\leq\ldots\leq\la_r$ such that the matrix of $J$ in $\B$ is given by{\footnotesize
\begin{equation}\label{diagonal}{\mathrm Mat}(J,\B)=\left(\begin{array}{cccccc}
\begin{array}{cc}0&-\la_1\\\la_1&0\end{array}&
0&\ldots&0&\ldots&0\\0&\ddots&\ddots&0&\ldots&0\\\ldots&\ddots&\ddots&0&\ldots&0\\
0&\ldots&0&\begin{array}{cc}0&-\la_r\\\la_r&0\end{array}&\ldots&0\\
\vdots&\vdots&\ldots&0&\ldots&0\\0&\ldots&\ldots&0&\ldots&0\end{array}\right).\end{equation}}
Let us return  to the general case.  Let $\B=(e_1,\bar e_2,\ldots,e_q,\bar e_q,f_1,\ldots,f_{n-2q})$ be a pseudo-Euclidean basis of $V$ and let $J\in Sym^-(V)$. Then its straightforward to see that the matrix of $J$ in $\B$ has the following form:
\begin{equation}\label{skew}\left\{\begin{array}{ccl}
{\mathrm Mat}(J,\B)&=&\left(\begin{array}{cc}A&P\\ \widehat{P} &B\end{array}\right),\;{}^tB=-B,\\
P&=&-\left(\begin{array}{c}X_1\\
Y_1\\
\vdots\\
X_q\\
Y_q\end{array}\right),\;\;
\widehat{P}=\left(\begin{array}{c}V_1\\
\vdots\\
V_{n-2q}\end{array}\right)=\left(\begin{array}{ccccc}{}^tY_1&{}^tX_1&
\ldots&{}^tY_q&{}^tX_q\end{array}\right),
\end{array}\right.\end{equation}where $X_i=(x_{1}^i,\ldots,x_{{n-2q}}^i)$,  $Y_i=(y_{1}^i,\ldots,y_{{n-2q}}^i)$, for $i=1,\ldots,q$ and $A=\left(A_{ij}\right)_{1\leq i,j\leq q}$ where the $A_{ij}$ are $(2,2)$-matrix  satisfying
\begin{equation}\label{amatrix}\left\{\begin{array}{ccl}A_{ii}&=&
\left(\begin{array}{cc}a_i&0\\0&-a_i\end{array}\right),\;i=1,\ldots,q,\\A_{ij}&=&
\left(\begin{array}{cc}a_{ij}&b_{ij}\\c_{ij}&d_{ij}\end{array}\right),\;
A_{ji}=
\left(\begin{array}{cc}-d_{ij}&-b_{ij}\\-c_{ij}&-a_{ij}\end{array}\right)\;\mbox{for}\; j>i.
\end{array}\right.\end{equation}

Let us give an expression of $\prs^*$ and compute its signature. If $J_1,J_2\in Sym^-(V)$  then
\begin{equation}\label{product}Mat(J_1\circ J_2,\B)=
\left(\begin{array}{cc}A_1A_2+P_1\widehat{P}_2&A_1P_2+P_1B_2\\\widehat{P}_1A_2+B_1
\widehat{P}_2&\widehat{P}_1P_2+B_1B_2
\end{array}\right).\end{equation}
One can see easily that
\begin{equation}\label{formule2}
P_1\widehat{P}_2=-\left(B_{ij}\right)_{1\leq i,j\leq q},\qquad
B_{ij}=\left(\begin{array}{cc}X^1_i.Y_j^2&X_i^1.X_j^2\\Y_i^1.Y_j^2&Y_i^1.X_j^2
\end{array}\right),
\end{equation}
\begin{equation}\label{formule3}
 \widehat{P}_1P_2=-\left(\langle V^1_i,V^2_j\rangle_q\right)_{1\leq i,j\leq n-2q},\end{equation}
where the dot is the canonical Euclidean product in $\R^{n-2q}$ and $\prs_q$ is the  pseudo-Euclidean product  of signature $(q,q)$ defined on $\R^{2q}$  by
    \begin{equation}\label{qproduct}\langle(x_1,y_1,\ldots,x_q,y_q),(x_1,y_1,\ldots,x_q,y_q)
    \rangle_q=2\sum_{i=1}^qx_iy_i.\end{equation}We shall denote by $\R^{(q,q)}$ the pseudo-Euclidean space $\R^{2q}$ endowed with $\prs_q$. From (\ref{product})-(\ref{formule3}) and the relation
\begin{equation}\label{dualformule}
\sum_{l=1}^q(X_l^1.Y_l^2+X_l^2.Y_l^1)=\sum_{l=1}^{n-2q}\langle V^1_l,V^2_l\rangle_q,\end{equation}one can deduce easily that
\begin{eqnarray}\label{trace}
{\mathrm tr}J_1\circ J_2&=&2\sum_{i=1}^qa_i^1a_i^2-2\sum_{l<k}(a_{lk}^1d_{lk}^2+d_{lk}^1a_{lk}^2+b_{lk}^1c_{lk}^2+
c_{lk}^1b_{lk}^2)\nonumber\\
&&-2\sum_{l=1}^q(X_l^1.Y_l^2+X_l^2.Y_l^1)+{\mathrm tr}B_1B_2.\end{eqnarray}

\begin{pr}\label{signature} The vector space $Sym^-(V)$  is of dimension $\frac{n(n-1)}2$ and the product $\prs^*$ is non degenerate and its  signature is $(q(n-q),\frac{n(n-1)+2q(q-n)}2)$.\end{pr}

 From what above, we can deduce that  $J\in Sym^-(V)$ is entirely determined, in  a pseudo-Euclidean basis,  by a $(n-2q,n-2q)$-matrix $B$ skew-symmetric in the Euclidean sense, a $(2q,2q)$-matrix $A$ skew-symmetric with respect to $\prs_q$ (which equivalent to $A$ satisfying (\ref{amatrix})), a family of vectors $(X_1,Y_1,\ldots,X_q,Y_q)$ in $\R^{n-2q}$ or a family of vectors $(V_1,\ldots,V_{n-2q})$ in $\R^{(q,q)}$.
We shall call, invariantly,  $$(A,B,X_1,Y_1,\ldots,X_q,Y_q)\quad\mbox{ or}\quad (A,B,V_1,\ldots,V_{n-2q})$$ the \emph{representation} of $J$ in $\B$. Note that we can choose the pseudo-Euclidean basis such that $B$ has the form (\ref{diagonal}).

\section{Ricci curvature of pseudo-Euclidean 2-step nilpotent Lie algebras}\label{section2}

A \emph{pseudo-Euclidean  Lie algebra}  is a pseudo-Euclidean vector space which is also a Lie algebra. Given a pseudo-Euclidean  Lie algebra $\G$, its associated \emph{Levi-Civita product}
is the bilinear map $\D:\G\times\G\too\G$ defined by the following relation:
\begin{equation}\label{levicivita}
2\langle\D_uv,w\rangle=\langle[u,v],w\rangle+\langle[w,u],v\rangle+\langle[w,v],u\rangle,\quad
u,v,w\in\G,\end{equation}
its associated curvature  is the three-linear map ${\cal R}:\G\times\G\times\G\too\G$ given by
$${\cal R}(u,v,w)=\D_{[u,v]}w-\D_u\D_vw+\D_v\D_uw,$$and the Ricci curvature is the bilinear map $\mathfrak{r}:\G\times\G\too\R$ given by
$$\mathfrak{r}(u,v)={\mathrm tr}(w\mapsto {\cal R}(u,w,v)).$$
 Denote by ${\cal J}:\G\too\G$ the endomorphism given by
$$\mathfrak{r}(u,v)=\langle {\cal J}\; u,v\rangle.$$ The scalar curvature of $\G$ is the real number
$\mathfrak{s}={\mathrm tr}{\cal J}.$\\
 A pseudo-Euclidean Lie algebra is called \emph{flat} (resp. \emph{Ricci flat}) if ${\cal R}=0$ (resp. ${\cal J}=0$).\\
A Lie algebra $\mathfrak{N}$ is called \emph{2-step nilpotent} if its derived ideal is non trivial and satisfies $[\mathfrak{N},\mathfrak{N}]\subset\mathfrak{Z}$, where $\mathfrak{Z}$ is the center of $\mathfrak{N}$. A trivial central extension of $\mathfrak{N}$ is a product of $\mathfrak{N}$ with an abelian Lie algebra. A 2-step nilpotent Lie algebra is called \emph{irreducible} if it is not a trivial central extension of any 2-step nilpotent Lie algebra. A Lie group is called 2-step nilpotent if its Lie algebra is a 2-step nilpotent Lie algebra.\\
Let $\mathfrak{N}$ be a pseudo-Euclidean 2-step nilpotent  Lie algebra and $(e_1,\ldots,e_p)$  a basis of
$\mathfrak{Z}$. Then, for any $u,v\in\mathfrak{N}$, the Lie bracket can be written
\begin{equation}\label{bracket}[u,v]=\sum_{i=1}^p\langle J_iu,v\rangle e_i,\end{equation}where $J_i:\mathfrak{N}\too\mathfrak{N}$ are skew-symmetric endomorphisms with respect to $\prs$ and
$\di\bigcap_{i=1}^p\ker J_i=\mathfrak{Z}.$ These endomorphisms will be called \emph{structure endomorphisms} associated to $(e_1,\ldots,e_p)$.
The structure endomorphisms $(K_1,\ldots,K_p)$ associated to a new basis $(f_1,\ldots,f_p)$  are given by
\begin{equation}\label{changementofbasis}
K_j=\sum_{i=1}^pp^{ji}J_i,\quad j=1,\ldots,p,\end{equation}where $(p^{ij})_{1\leq i,j\leq p}$ is  the passage matrix from $(f_1,\ldots,f_p)$ to $(e_1,\ldots,e_p)$.

The proof of the following proposition is  a direct computation.
\begin{pr}Let $\mathfrak{N}$ be a pseudo-Euclidean 2-step nilpotent  Lie algebra, $(e_1,\ldots,e_p)$  a basis of
$\mathfrak{Z}$ and $(J_1,\ldots,J_p)$ the corresponding structure endomorphisms. Then the endomorphisms $F,G:\mathfrak{N}\too\mathfrak{N}$  given by
\begin{equation}\label{invariant}
F=\frac12\sum_{i,j=1}^p\langle e_i,e_j\rangle J_i\circ J_j\quad
\mbox{and}\quad G(u)=-\frac14\sum_{i,j=1}^p\langle e_i,u\rangle{\mathrm tr}(J_i\circ J_j)e_j\end{equation}are symmetric with respect to $\prs$, are independent of the choice of the basis $(e_1,\ldots,e_p)$ and satisfy $F\circ G=G\circ F=0$. \end{pr}

 {\bf Notation.} We shall denote by ${\cal J}^- $ and ${\cal J}^+$, respectively, the symmetric endomorphisms $F$ and $G$ defined by (\ref{invariant}). We denote also by $\mathfrak{r}^-$ and $\mathfrak{r}^+$ the symmetric bilinear forms
$$\mathfrak{r}^-(u,v)=\langle {\cal J}^-u,v\rangle \quad
\mbox{and}\quad \mathfrak{r}^+(u,v)=\langle {\cal J}^+u,v\rangle.$$
These notations are justified by  the following lemma which  will play a crucial role in this paper.
\begin{Le}\label{mainlemma}Let $\mathfrak{N}$ be a pseudo-Euclidean 2-step nilpotent  Lie algebra. Then its Ricci curvature  is given by
\begin{equation}\label{ricci}\mathfrak{r}=\mathfrak{r}^++\mathfrak{r}^-,
\end{equation}and its scalar curvature  is given by
\begin{equation}\label{scalar}\mathfrak{s}=\frac12{\mathrm tr}{\cal J}^-=-{\mathrm tr}{\cal J}^+.\end{equation}\end{Le}
{\it Proof.} Note first that if $(J_1,\ldots,J_p)$ are the structure endomorphisms associated to a  basis  $(e_1,\ldots,e_p)$  of
$\mathfrak{Z}$, one can deduce easily form (\ref{levicivita}) and (\ref{bracket}) that the Levi-Civita product of $\prs$ is given by
\begin{eqnarray*}
2\D_uv&=&\sum_{i=1}^p\left(\langle J_iu,v\rangle  e_i
-  \langle e_i,v\rangle J_iu
-  \langle e_i,u\rangle J_iv\right),\end{eqnarray*}and its curvature is given by
{\small\begin{eqnarray}
R(u,v)w&=&\sum_{i,j=1}^p\langle e_i,e_j\rangle\left(\frac14\langle J_iv,w\rangle J_ju -\frac14\langle J_iu,w\rangle   J_jv-\frac12\langle J_iu,v\rangle  J_jw\right)\nonumber\\
&&+\frac14\sum_{i,j=1}^p\left(\langle e_i,w\rangle\langle e_j,v\rangle J_j\circ J_iu-
\langle e_i,w\rangle\langle e_j,u\rangle J_j\circ J_iv+
\langle e_i,u\rangle
\langle e_j,v\rangle [J_j, J_i]w\right)\nonumber\\
&&+\frac14\sum_{i,j=1}^p\left(\langle e_i,w\rangle \langle [J_j,J_i]u,v\rangle
+\langle e_i,v\rangle \langle J_ju,J_iw\rangle-  \langle e_i,u\rangle \langle J_jv,J_iw\rangle\right)e_j.\label{curvaturenilpotent}
\end{eqnarray}}
 Let $\mathfrak{F}$ (resp. $\mathfrak{G}$) be a complement of $\mathfrak{Z}\cap\mathfrak{Z}^\perp$ in $\mathfrak{Z}$ (resp. in $\mathfrak{Z}^\perp$). The subspace  $\mathfrak{F}\oplus \mathfrak{G}$ is nondegenerate and its orthogonal, which is also nondegenerate, has dimension $2\dim\mathfrak{Z}\cap\mathfrak{Z}^\perp$. So we can construct a basis $$(e_1,\ldots,e_q,\bar e_1,\ldots,\bar e_q,f_1,\ldots,f_{p-q},g_1,\ldots,g_{n-p-q})$$ of $\mathfrak{N}$ such that:
\begin{enumerate}\item $(e_1,\ldots,e_q,\bar e_1,\ldots,\bar e_q)$ is a basis of $\mathfrak{F}^\perp\cap \mathfrak{G}^\perp$ and, for  $i,j=1,\ldots,q,$ $\langle e_i,\bar e_j\rangle=\de_{ij},$
$\langle e_i, e_j\rangle=0$ and $\langle \bar e_i,\bar e_j\rangle=0,$
\item $(f_1,\ldots,f_{p-q})$ and  $(g_1,\ldots,g_{n-p-q})$ are orthogonal basis, respectively, of $\mathfrak{F}$ and $\mathfrak{G}$, moreover, for $i=1,\ldots,p-q$ and $j=1,\ldots,n-p-q$, $\langle f_i,f_i\rangle=\pm1$ and $\langle g_j,g_j\rangle=\pm1$.

\end{enumerate}
The Ricci curvature is given by
\begin{eqnarray*}
\mathfrak{r}(u,v)&=&\sum_{i=1}^q\left(\langle R(u,e_i)v,\bar e_i\rangle+
\langle R(u,\bar e_i)v, e_i\rangle\right)+\sum_{i=1}^{p-q}\langle R(u,f_i)v, f_i\rangle\langle f_i,f_i\rangle\\&&+\sum_{i=1}^{n-p-q}\langle R(u,g_i)v, g_i\rangle\langle g_i,g_i\rangle.\end{eqnarray*}

We denote by $(K_1,\ldots,K_q,J_1,\ldots,J_{p-q})$ the structure endomorphisms associated to $(e_1,\ldots,e_{q},f_1,\ldots,f_{p-q})$.
By using (\ref{curvaturenilpotent}), we get
\begin{eqnarray*}
R(u,e_i)v&=&0,\\
\langle R(u,\bar e_i)v, e_i\rangle&=&0,\\
\langle R(u,f_k)v, f_k\rangle&=&\frac14
 \langle J_ku,J_kv\rangle,\\
 \langle R(u,g_k)v, g_k\rangle&=&-\frac34\sum_{i=1}^{p-q}\langle f_i,f_i\rangle\langle g_k,J_iv\rangle \langle J_iu, g_k\rangle\\
 &&-\frac14\sum_{i,j=1}^q
\langle e_i,v\rangle\langle e_j,u\rangle \langle K_j\circ K_ig_k,g_k\rangle
-\frac14\sum_{i,j=1}^{p-q}
\langle f_i,v\rangle\langle f_j,u\rangle \langle J_j\circ J_ig_k,g_k\rangle\\
&&-\frac14\sum_{i,j=1}^{q,p-q}
\langle e_i,v\rangle\langle f_j,u\rangle \langle J_j\circ K_ig_k,g_k\rangle
-\frac14\sum_{j,i=1}^{q,p-q}\langle e_j,v\rangle\langle f_i,u\rangle \langle J_i\circ K_jg_k,g_k\rangle,
\end{eqnarray*}and (\ref{ricci}) follows. Let us establish (\ref{scalar}). We have, from (\ref{ricci}),
$$\mathfrak{s}={\mathrm tr}{\cal J}^++{\mathrm tr}{\cal J}^-.$$ On the other hand,
\begin{eqnarray*}
{\mathrm tr}{\cal J}^-&\stackrel{(\ref{invariant})}=&\frac12\sum_{i=1}^{p-q}\langle f_i,f_i\rangle {\mathrm tr}(J_i^2),\\
{\mathrm tr}{\cal J}^+&=&2\sum_{i=1}^q\langle {\cal J}^+e_i,\bar e_i\rangle+
\sum_{i=1}^{p-q}\langle {\cal J}^+f_i,f_i\rangle\langle f_i,f_i\rangle
+\sum_{i=1}^{n-p-q}\langle {\cal J}^+g_i,g_i\rangle\langle g_i,g_i\rangle\\
&\stackrel{(\ref{invariant})}=&-\frac14\sum_{i=1}^{p-q}{\mathrm tr}(J_i^2)\langle f_i,f_i\rangle.\end{eqnarray*}So we deduce
${\mathrm tr}{\cal J}^+=-\frac12{\mathrm tr}{\cal J}^-$ and (\ref{scalar}) follows.
\hfill $\;\square$

\begin{rem}\begin{enumerate}\item The Ricci curvature and the scalar curvature of a pseudo-Euclidean 2-step nilpotent Lie algebra were computed in \cite{cordero} (see Theorems 3.24 and 3.26). The formulas (\ref{ricci}) and (\ref{scalar}) are more simple than those given in \cite{cordero}.
\item The situation is quite simple in the Euclidean case. Indeed, for any orthonormal basis $(e_1,\ldots,e_p)$ of $\mathfrak{Z}$, and any orthonormal basis $(f_1,\ldots,f_{n-p})$ of $\mathfrak{Z}^\perp$, we have
    $${\cal J}^-=\frac12\sum_{j=1}^pJ_j^2$$ {and} the matrix of ${\cal J}^+$ in $(e_1,\ldots,e_p,f_1,\ldots,f_{n-p})$ is given by $$-\frac14\left(\begin{array}{cc}\left({\mathrm tr}(J_i\circ J_j)\right)_{1\leq i,j\leq p}&0\\0&0\end{array}\right),$$ where $(J_1,\ldots,J_p)$ are the structure endomorphisms associated to $(e_1,\ldots,e_p)$.
So the notations $\mathfrak{r}^+$ and $\mathfrak{r}^-$ in (\ref{ricci}) are appropriate in  this  case.  However, in the general case, ${\cal J}^+$ (resp.${\cal J}^-$) can have negative (resp. positive) eigenvalues or non real eigenvalues. Nevertheless, we will conserve these notations.
\end{enumerate}\end{rem}

The following proposition is a more accurate version of a result of Eberlein (see \cite{Eberlein} Proposition 2.5) and can be deduced easily from what above.
\begin{pr}Let $\mathfrak{N}$ be an Euclidean  2-step nilpotent Lie algebra of dimension $n$. Put $p=\dim\mathfrak{Z}$ and $r=\dim[\mathfrak{N},\mathfrak{N}]$. Then:
 \begin{enumerate}\item ${\mathrm rank}{\cal J}^+=r$ and ${\mathrm rank}{\cal J}^-=n-p$,
 \item there exists an orthonormal basis $(e_1,\ldots,e_p)$ of $\mathfrak{Z}$,  an orthonormal  basis $(g_1,\ldots,g_{n-p})$ of $\mathfrak{Z}^\perp$ and two families of real numbers $0<\mu_1\leq \ldots\leq \mu_{r}$ and $0<\la_1\leq \ldots\leq \la_{n-p}$ such that non vanishing entries in the  matrix of the Ricci curvature $\mathfrak{r}$ in the basis $\B=(e_1,\ldots,e_p,g_1,\ldots,g_{n-p})$ are
     $$\mathfrak{r}(e_i,e_i)=\mu_i \quad\mbox{and}\quad
     \mathfrak{r}(g_j,g_j)=-\la_j,\; i=1,\ldots,r,\;j=1,\ldots,n-p.$$
Moreover, the scalar curvature is given by
\[\label{scalareucldian}
\mathfrak{s}=-\frac12(\la_1+\ldots+\la_{n-p})=-(\mu_1+\ldots+\mu_r).\]In particular, the scalar curvature is negative.
\end{enumerate}
\end{pr}
The scalar curvature of an Euclidean 2-step nilpotent Lie algebra is negative which is expectable according to a theorem by Uesu (see \cite{uesu}).

There is an important class of Euclidean 2-step nilpotent Lie algebras introduced by Kaplan in \cite{kaplan} and called \emph{Heisenberg type}  Lie algebras. Let us compute the Ricci curvature and the scalar curvature of such Lie algebras. Recall that a {Heisenberg type}  Lie algebra is a  2-step nilpotent Lie algebra $\mathfrak{N}$ endowed with an Euclidean product $\prs$ such that, for any $z\in\mathfrak{Z}$ and for any $v\in\mathfrak{Z}^\perp$,
\begin{equation}\label{htype}\langle ad_{v}^tz,ad_{v}^tz\rangle=\langle z,z\rangle\langle v,v\rangle,\end{equation}where $ad_u^t$ is the adjoint of $ad_u$ with respect to $\prs$.
\begin{pr}Let $\mathfrak{N}$ be a Heisenberg type Lie algebra. Then
\[    {\cal J}^+=\frac14\dim\mathfrak{Z}^\perp P_{\mathfrak{Z}},\quad
{\cal J}^-=-\frac12\dim\mathfrak{Z}P_{\mathfrak{Z}^\perp}\quad\mbox{
and}\quad
\mathfrak{s}=-\frac14\dim\mathfrak{Z}\dim\mathfrak{Z}^\perp,\]where
$P_{\mathfrak{Z}^\perp}$ and $P_{\mathfrak{Z}}$ denote the orthogonal projections onto $\mathfrak{Z}^\perp$ and $\mathfrak{Z}$ respectively.

\end{pr}
{\it Proof.} Choose an orthonormal basis $(e_1,\ldots,e_p)$ of $\mathfrak{Z}$ and denote by $(J_1,\ldots,J_p)$ the associated structure endomorphisms. One can see easily that, for any $z\in\mathfrak{Z}$ and for any $v\in\mathfrak{Z}^\perp$,
$$ad_{v}^tz=\sum_{i=1}^p\langle z,e_i\rangle J_iv.$$Thus (\ref{htype}) is equivalent to
\[\sum_{k,l=1}^p\langle e_k,e_i\rangle\langle e_l,e_j\rangle
\langle J_kv,J_lv\rangle=\langle e_i,e_j\rangle\langle v,v\rangle,\quad i,j=1,\ldots,p, v\in\mathfrak{Z}^\perp,\]which is equivalent to
\[
J_i^2=-P_{\mathfrak{Z}^\perp},\; i=1,\ldots,p\quad\mbox{and}\quad J_i\circ J_j=-J_j\circ J_i, i\not=j.
\]
The proposition follows from these relations and (\ref{invariant})-(\ref{scalar}). \hfill $\;\square$

Let us show now that an Einstein pseudo-Euclidean 2-step nilpotent Lie algebra must be Ricci flat.
\begin{pr} Let $\mathfrak{N}$ be a pseudo-Euclidean 2-step nilpotent Lie algebra such that there exists $\la\in\R$
satisfying $ \mathfrak{r}=\la\prs$. Then $\la=0$.
\end{pr}

{\it Proof.} We consider the basis $(e_1,\ldots,e_q,\bar e_1,\ldots,\bar e_q,f_1,\ldots,f_{p-q},g_1,\ldots,g_{n-p-q})$ of $\mathfrak{N}$ constructed in the proof of Lemma \ref{mainlemma} and we denote by $(K_1,\ldots,K_q,J_1,\ldots,J_{p-q})$ the structure endomorphisms associated to $(e_1,\ldots,e_{q},f_1,\ldots,f_{p-q})$.

We distinguish two cases:
\begin{enumerate}\item The center is degenerate. In this case $0=({\cal J}^-+{\cal J}^+)e_i=\la e_i$ and hence $\la=0$.
\item The center is nondegenerate. By using (\ref{invariant}) and (\ref{ricci}), we can see that $\mathfrak{r}=\la \prs$ is equivalent to
    $$
    \left\{\begin{array}{c}
    \di\frac12\sum_{i=1}^{p}\langle f_i,f_i\rangle J_i^2g_j=\la g_j,\quad j=1,\ldots,n-p-q,\\
    -\frac14\langle f_i,f_i\rangle{\mathrm tr}(J_i^2)f_i=\la f_i,\quad i=1,\ldots,p-q.\end{array}\right.
    $$
    Thus, for $i=1,\ldots,p-q$ and $j=1,\ldots,n-p-q$,
    $$\la=-\frac14\langle f_i,f_i\rangle{\mathrm tr}(J_i^2)=
    \frac12\sum_{i=1}^{p}\langle f_i,f_i\rangle \langle g_j,g_j\rangle \langle J_i^2g_j,g_j\rangle$$ and hence
    $$\dim\mathfrak{Z}^\perp\la=\frac12\sum_{i=1}^{p}\langle f_i,f_i\rangle
    {\mathrm tr}(J_i^2)=-2\dim\mathfrak{Z}\la$$which implies $\la=0$.\hfill $\square$

\end{enumerate}

\section{ Ricci flat left invariant pseudo-Riemannian metrics on Heisenberg groups}\label{section3}

A Heisenberg Lie algebra is 2-step nilpotent Lie algebra   of dimension $2k+1$ such that its center is of dimension 1 and coincides with its derived ideal. A $2k+1$-dimensional Heisenberg Lie algebra is isomorphic to
$${\h}_{2k+1}=\left\{\left(\begin{array}{ccc}0&X&z\\0&0&{}^t\bar{X}\\0&0&0\end{array}\right),z\in\R,X,\bar{X}\in
\R^k\right\},$$whose associated simply connected Lie group is
$$H_{2k+1}=\left\{\left(\begin{array}{ccc}1&X&z\\0&I&{}^t\bar{X}\\0&0&1\end{array}\right),z\in\R,X,\bar{X}\in
\R^k\right\}.$$We shall denote by $(z,x_1,\bar{x}_1,\ldots,x_k,\bar{x}_k)$ the canonical coordinates of $H_{2k+1}$.

It is well-known that a $2k+1$-dimensional Heisenberg Lie group carries a  flat  pseudo-Riemannian  metric if and only if $k=1$ and in this case the metric is Lorentzian (see \cite{medina}).
In this section, we resolve completely the problem of existence of Ricci flat pseudo-Riemannian  metrics on Heisenberg Lie groups. We  show the following result.
\begin{theo}\label{heisenberg}Let $\G$ be a $2k+1$-dimensional Heisenberg Lie algebra. Then:\\$(i)$ $\G$ carries a
 Ricci flat  Lorentz product if and only if $k=1$.\\
$(ii)$ For any $q$ such that
$2\leq q\leq k$, $\G$ carries a Ricci flat  pseudo-Euclidean product of signature $(q,2k+1-q)$.

\end{theo}

 {\it Proof.} First, if $\prs$ is a pseudo-Euclidean  product on $\G$, $e$ a non null central element  and  $J$ the associated structure endomorphism, then
$${\cal J}^+=-\frac14\langle e,\;.\;\rangle {\mathrm tr}(J^2)e\quad\mbox{and}\quad
{\cal J}^-=\frac12\langle e,e\rangle J^2,$$
and one can see easily from (\ref{ricci}) that $\mathfrak{r}=0$ if and only if
\begin{equation}\label{ricciheisenberg}\langle e,e\rangle={\mathrm tr}(J^2)=0.\end{equation}

$(i)$
We suppose that $\G$ carries a  Ricci flat Lorentzian product, we choose a non null central element $e$ and
we denote by $J$ the associated structure endomorphism. Let $\B=(e,\bar e,f_1,\ldots,f_{2k-1})$ be a Lorentzian
basis of $\G$. Since $\ker J=\R e$, we deduce from (\ref{skew}) that the representation of $J$ in $\B$ is
$(0,B,X_1,0)$, where $B$ is a skew-symmetric $(2k-1,2k-1)$-matrix and $X_1=(x_1,\ldots,x_{2k-1})$. On the other hand,
we deduce from (\ref{trace}) that ${\mathrm tr}(J^2)={\mathrm tr}(B^2)$ and hence (\ref{ricciheisenberg}) implies $B=0$. Thus the representation of $J$ in $\B$ is $(0,0,X_1,0)$ and the condition $\ker J=\R e$ implies $2k-1=1$, thus $k=1$. Conversely, if $k=1$, it is well-known that $\h_3$ carries a  flat Lorentzian  product.

$(ii)$ Suppose that $2\leq q\leq k$.
We consider  the pseudo-Euclidean space
$V=\R^{2q}\times\R^{2(k-q)+1}$ endowed with the inner product of signature $(q,2k+1-q)$ given by
$$\langle(u,v),(u,v)\rangle=\langle u,u\rangle_q+v.v,$$where $\prs_q$ is the  pseudo-Euclidean product on $\R^{2q}$ given by (\ref{qproduct}) and the dot is the canonical Euclidean product on $\R^{2(k-q)+1}$. For any $U=(x_1,y_1,\ldots,x_q,y_q)\in\R^{2q}$, we put
    $\overline{U}=(y_1,x_1,\ldots,y_q,x_q)$ and we denote by $(e_1,\bar e_1,e_2,\bar e_2,\ldots,e_q,\bar e_q)$  the canonical basis of $\R^{2q}$.

We will construct a skew-symmetric endomorphism $J$ on $V$ such that ${\mathrm tr}(J^2)=0$ and $\ker J=\R e_1$.  Once $J$ is constructed, the Lie bracket $[u,v]=\langle Ju,v\rangle e_1$ shall induce on $V$ a structure of Heisenberg's Lie algebra  for which the pseudo-Euclidean product $\prs$ is Ricci flat. We  give a general method for building such an endomorphism. Explicit examples will be given in Example \ref{example1}.

We distinguish  two cases:
    \begin{enumerate}\item $k=2q+r$ with $r\geq0$. Let $J$ be the skew-symmetric endomorphism of $V$ whose matrix in the canonical basis is $\left(\begin{array}{cc}A&P\\\widehat{P}&B\end{array}\right)$ where:
    \begin{enumerate}\item $\widehat{P}=\left(\begin{array}{ccccccc}0&{}^tX_1&{}^tY_2&^tX_2&\ldots&{}^tY_q&{}^tX_q\end{array}\right)$,\item $(X_1,X_2,Y_2,\ldots,X_q,Y_q)$ is a family of linearly independent vectors in $\R^{2(k-q)+1}=\R^{2q+2r+1}$,
    \item $\mbox{rank}B=2(r+1)$ and
    \begin{equation}\label{first}span\{X_1,X_2,Y_2,\ldots,X_q,Y_q\}\oplus{\mathrm Im}B=\R^{2(k-q)+1},\end{equation}
    \item $Ae_1=0$ and
    \begin{equation}\label{second}{\mathrm tr}(A^2)+{\mathrm tr}(B^2)=4\sum_{i=2}^qX_i.Y_i.\end{equation}

    \end{enumerate}
    One can deduce from (\ref{trace}) and (\ref{second}) that ${\mathrm tr}(J^2)=0$. Let us check now that $\ker J=\R e_1$.
     A vector $(a_1,b_1,\ldots,a_q,b_q,Z)\in\ker J$ if and only if
    $$\left\{\begin{array}{ccccc}AU&=&\left( X_1.Z, Y_1.Z,\ldots,
 X_q.Z, Y_q.Z\right)
&&\\BZ&=&\di-b_1X_1-\sum_{i=2}^q(a_iY_i+b_iX_i)&&\end{array}\right.$$
where $U=(a_1,b_1,\ldots,a_q,b_q)$. From (\ref{first}) one can deduce that $U=(a_1,0,\ldots,0)$ and $Z\in\ker B\cap span\{X_1,X_2,Y_2,\ldots,X_q,Y_q\}^\perp=\{0\}$ and the result follows.

    \item $k=q+r$ with $0\leq r\leq q-1$. Let $J$ be the skew-symmetric endomorphism of $V$ whose matrix in the canonical basis is $\left(\begin{array}{cc}A&P\\\widehat{P}&B\end{array}\right)$ where:
        \begin{enumerate} \item $\widehat{P}=\left(\begin{array}{c}\overline{V}_1\\\vdots\\\overline{V}_{2r+1}\end{array}\right)$,
        \item $(V_1,\ldots,V_{2r+1})$ is a family of linearly independent vectors in $$F=Vect\{e_1,e_2,\bar e_2,\ldots,e_q,\bar e_q\}\subset\R^{2q},$$
            \item $\mbox{rank}A=2(q-r-1)$ and
    \begin{eqnarray}\label{first1}span\{V_1,\ldots,V_{2r+1}\}\oplus {\mathrm Im}A=F,\\
    \label{second1}{\mathrm tr}(A^2)+{\mathrm tr}(B^2)=2\sum_{i=1}^{2r+1}\langle V_i,V_i\rangle_q.
    \end{eqnarray}\end{enumerate}
    One can deduce from (\ref{trace}) and (\ref{second1}) that ${\mathrm tr}(J^2)=0$. Let us check now that $\ker J=\R e_1$.
    A vector $(U,Z)=(a_1,b_1,\ldots,a_q,b_q,z_1,\ldots,z_{2r+1})\in\ker J$ if and only if
    $$\left\{\begin{array}{ccccc}AU&=&\di\sum_{i=1}^{2r+1}z_iV_i,
&&\\BZ&=&-\left(\langle V_1,U\rangle_q,\ldots,
\langle V_{2r+1},U\rangle_q\right)&&\end{array}\right.$$From (\ref{first1}), we deduce that $Z=0$ and hence $U\in\ker A\cap span\{V_1,\ldots,V_{2r+1}\}^\perp=F^\perp=\R e_1$ and the result follows.
\hfill $\;\square$

\end{enumerate}

\begin{exem}\label{example1}We give three explicit examples illustrating  the general construction done in the proof of Theorem \ref{heisenberg} and, moreover, for each case we give the expression of the corresponding left invariant pseudo-Riemannian metric on the associated simply connected Lie group.
\begin{enumerate}\item  $2\leq q$ and $k=2q+r$. We consider $V=\R^{2q}\times\R^{2q-1}\times\R^{2(r+1)}$ endowed with the product
$$\langle(u,v,w),(u,v,w)\rangle=\langle u,u\rangle_q+v.v+w.w.$$We denote by $(e_1,\bar{e}_1,\ldots,e_q,\bar{e}_q)$, $(f,f_1,\bar{f}_1,\ldots,f_{q-1},\bar{f}_{q-1})$ and $(g_1,\bar{g}_1,\ldots,g_{r+1},\bar{g}_{r+1})$ the canonical basis of $\R^{2q}$, $\R^{2q-1}$ and $\R^{2(r+1)}$, respectively.
We consider the skew-symmetric endomorphism $J$ of $V$ given by
 \begin{eqnarray*}
 Je_1&=&0,\;
 J\bar{e}_1=-f,\;
 Je_i=a_ie_i-\bar{f}_{i-1},\;
 J\bar{e}_i=-a_i\bar{e}_i-f_{i-1}, i=2,\ldots,q,\\
 Jf&=&e_1,\;
 Jf_{i}=e_{i+1},\; J\bar{f}_i=\bar{e}_{i+1},\; i=1,\ldots,q-1,\\
 Jg_i&=&\la_i \bar{g}_i,\; J\bar{g}_i=-\la_ig_i,\; i=1,\ldots,r+1,\end{eqnarray*}where $0<\la_1\leq\ldots\leq\la_{r+1}$ and $\di\sum_{i=2}^qa_i^2=\sum_{i=1}^{r+1}\la_i^2.$ One can check easily that $\ker J=\R e_1$ and ${\mathrm tr}(J^2)=0$. Thus  the bracket $[u,v]=\langle Ju,v\rangle e_1$ induces on $V$ a Heisenberg's Lie algebra structure for which $\prs$ is Ricci flat.
 Put
 \begin{eqnarray*}
 E&=&e_1,\; F_1=f,\; \bar{F}_1=\bar{e}_1,\;{F}_i=f_{i-1},\;\bar F_i=\bar{e}_i+a_i\bar{f}_{i-1},\; i=2,\ldots,q,\\
 E_i&=&\bar{f}_{i},\;\bar{E}_i=e_{i+1},\; i=1,\ldots,q-1,\\
 G_i&=&\frac1{\sqrt{\la_i}}g_i,\;\bar{G}_i=\frac1{\sqrt{\la_i}}\bar{g}_i,\; i=1,\ldots,r+1.\end{eqnarray*}
$(E,E_1,\bar{E}_1,\ldots,E_{q-1},\bar{E}_{q-1},F_1,\bar{F}_1,\ldots,F_q,\bar{F}_q,G_1,\bar{G}_1,\ldots,G_{r+1},\bar{G}_{r+1})$
is a basis of $V$ satisfying
$$\label{bracketexam}[E_i,\bar{E}_i]=[F_j,\bar{F}_j]=[G_k,\bar{G}_k]=E,\;i=1,\ldots,q-1,\; j=1,\ldots,q\quad\mbox{and}\quad k=1,\ldots,r+1,$$ and the all the others brackets are null or given by anti-symmetry. Thanks to this basis we identify $(V,[\;,\;],\prs)$ with the Lie algebra
$${\h}_{2(2q+r)+1}=\left\{\left(\begin{array}{ccc}0&X&z\\0&0&{}^t\bar{X}\\0&0&0\end{array}\right),z\in\R,X,\bar{X}\in\R^{q-1}
\times\R^{q}\times\R^{r+1}\right\},$$endowed with the pseudo-Euclidean product of signature $(q,2(2q+r)+1-q)$ whose expression in the canonical basis $(e,\ell_i,\bar{\ell}_i,m_j,\bar{m}_j,n_k,\bar{n}_k)$ $i=1,\ldots,q-1$, $j=1,\ldots,q$, $k=1,\ldots,r+1$ of
${\h}_{2(2q+r)+1}$ is given by
\begin{eqnarray*}
\langle m_i,m_i\rangle&=&\langle \ell_j,\ell_j\rangle=
\langle e,\bar{m}_1\rangle=1,\; i=1,\ldots,q,\; j=1,\ldots,q-1,\\
\langle \bar \ell_{j-1},\bar m_j\rangle&=&1,\; j=2,\ldots,q,\\
\langle n_i,n_i\rangle&=&\langle \bar n_i,\bar n_i\rangle=\la_i^{-1},\; i=1,\ldots,r+1,\\
\langle \bar m_i,\bar m_i\rangle&=&a_i^2,\;\langle  \ell_{i-1},\bar m_i\rangle=a_i,\; i=2,\ldots,q,\end{eqnarray*}and all the others products are null.
The simply connected Lie group associated to ${\h}_{2(2q+r)+1}$ is given by
$${H}_{2(2q+r)+1}=\left\{\left(\begin{array}{ccc}1&X&z\\0&I&{}^t\bar{X}\\0&0&1\end{array}\right),z\in\R,X,\bar{X}\in\R^{q-1}
\times\R^{q}\times\R^{r+1}\right\},$$and if
$(z,x_i,\bar{x}_i,y_j,\bar{y}_j,t_k,\bar{t}_k)$ are its canonical coordinates then the left invariant pseudo-Riemannian metric $\prs^l$ on ${H}_{2(2q+r)+1}$ associated to $\prs$ is given by
\begin{eqnarray*}
\prs^l&=&2d\bar{y_1}(dz-\sum_{i=1}^{q-1}x_id\bar{x_i}-\sum_{i=1}^qy_id\bar{y}_i-
\sum_{i=1}^{r+1}t_id\bar{t}_i)+\sum_{i=1}^qdy_i^2\\&&+\sum_{i=2}^q\left(2d\bar{x}_{i-1}d\bar{y}_i+(dx_{i-1}+a_id\bar{y}_i)^2\right)+
\sum_{i=1}^{r+1}\la_i^{-1}(dt_i^2+d\bar{t}_i^2),\end{eqnarray*}
where $0<\la_1\leq\ldots\leq\la_{r+1}$ and $\di\sum_{i=2}^qa_i^2=\sum_{i=1}^{r+1}\la_i^2.$
\item $2\leq q$ and $k=q+r$ with $1\leq r\leq q-1$. Let $V=\R^{2q}\times\R^{2r+1}$ endowed with the product
$$\langle(u,v),(u,v)\rangle=\langle u,u\rangle_q+v.v.$$We denote by $(e_1,\bar{e}_1,\ldots,e_q,\bar{e}_q)$ and $(f,f_1,\bar{f}_1,\ldots,f_{r},\bar{f}_{r})$  the canonical basis of $\R^{2q}$ and $\R^{2r+1}$, respectively.
We consider the skew-symmetric endomorphism $J$ of $V$ given by
 \begin{eqnarray*}
 Je_1&=&0,\;
 J\bar{e}_1=-f,\;\\
 Je_i&=&-\bar{f}_{i-1},\;
 J\bar{e}_i=-f_{i-1}, i=2,\ldots,r+1,\\
 Je_j&=&a_je_j,\;
 J\bar{e}_j=-a_j\bar{e}_j, j=r+2,\ldots,q,\\
 Jf&=&e_1,\;
 Jf_{i}=e_{i+1}+\la_i\bar{f}_i,\; J\bar{f}_i=\bar{e}_{i+1}-\la_if_i,\; i=1,\ldots,r,
 \end{eqnarray*}where  $0<a_{r+2}\leq\ldots\leq a_q$ and $\di\sum_{i=r+2}^qa_i^2=\sum_{i=1}^{r}\la_i^2$ if $r<q-1$ and $\la_1=\ldots=\la_r=0$ if $r=q-1$. One can check easily that $\ker J=\R e_1$ and ${\mathrm tr}(J^2)=0$. Thus  the bracket $[u,v]=\langle Ju,v\rangle e_1$ induces on $V$ a Heisenberg's Lie algebra structure for which $\prs$ is Ricci flat.
 Put
 \begin{eqnarray*}
 E&=&e_1,\; F_1=f,\; \bar{F}_1=\bar{e}_1,\;{F}_i=\la_i e_i+f_{i-1},\;\bar F_i=\bar{e}_i,\; i=2,\ldots,r+1\\
 F_i&=&\frac1{\sqrt{a_i}}e_i,\;\bar{F}_i=\frac1{\sqrt{a_i}}\bar{e}_i,\; i=r+2,\ldots,q,\\
 E_i&=&\bar{f}_{i},\;\bar{E}_i=e_{i+1},\; i=1,\ldots,r.
 \end{eqnarray*}
 $(E,E_1,\bar{E}_1,\ldots,E_{r},\bar{E}_{r},F_1,\bar{F}_1,\ldots,F_q,\bar{F}_q)$
is a basis of $V$ satisfying
$$\label{bracketexam}[E_i,\bar{E}_i]=[F_j,\bar{F}_j]=E,\;i=1,\ldots,r,\; j=1,\ldots,q,$$ and the all the others brackets are null or given by anti-symmetry. Thanks to this basis we identify $(V,[\;,\;],\prs)$ with the Lie algebra
$${\h}_{2(q+r)+1}=\left\{\left(\begin{array}{ccc}0&X&z\\0&0&{}^t\bar{X}\\0&0&0\end{array}\right),z\in\R,X,\bar{X}\in\R^{r}
\times\R^{q}\right\},$$endowed with the pseudo-Euclidean product of signature $(q,2(2q+r)+1-q)$ whose expression in the canonical basis $(e,\ell_i,\bar{\ell}_i,m_j,\bar{m}_j)$ $i=1,\ldots,r$, $j=1,\ldots,q$ of
${\h}_{2(q+r)+1}$ is given by
 \begin{eqnarray*}
 \langle e,\bar{m}_1\rangle&=&\langle \ell_i,{\ell}_i\rangle=\langle m_j,{m}_j\rangle=1,\; i=1,\ldots,r,\;j=1,\ldots,r+1,\\
 \langle \bar{m}_j,\bar{\ell}_{j-1}\rangle&=&1,\; j=2,\ldots,r+1,\\
 \langle m_i,\bar{m}_i\rangle&=&\la_i,\; i=2,\ldots,r+1,\\
 \langle m_i,\bar{m}_i\rangle&=&a_i^{-1},\;i=r+2,\ldots,q,\end{eqnarray*}and all the others products are null.
 The simply connected Lie group associated to ${\h}_{2(q+r)+1}$ is given by
$${H}_{2(q+r)+1}=\left\{\left(\begin{array}{ccc}1&X&z\\0&I&{}^t\bar{X}\\0&0&1\end{array}\right),z\in\R,X,\bar{X}\in\R^{r}
\times\R^{q}\right\},$$and if
$(z,x_i,\bar{x}_i,y_j,\bar{y}_j)$ are its canonical coordinates then the left invariant pseudo-Riemannian metric $\prs^l$ on ${H}_{2(q+r)+1}$ associated to $\prs$ is given by
\begin{eqnarray*}
 \prs^l&=&2d\bar{y}_1\left(dz-\sum_{i=1}^{r}x_id\bar{x}_i-\sum_{i=1}^{q}y_id\bar{y}_i\right)
 +\sum_{i=1}^{r+1}dy_i^2+\sum_{i=1}^{r}\left(dx_i^2+2d\bar{x_i}d\bar{y}_{i+1}\right)\\&&
 +2\sum_{i=1}^{r+1}\la_idy_id\bar{y_i}+2\sum_{i=r+2}^{q}a_i^{-1}dy_id\bar{y_i},\end{eqnarray*}where  $0<a_{r+2}\leq\ldots\leq a_q$ and $\di\sum_{i=r+2}^qa_i^2=\sum_{i=1}^{r}\la_i^2$ if $r<q-1$ and $\la_1=\ldots=\la_r=0$ if $r=q-1$.
\item $ k=q$ and $q\geq2$. Let $V=\R^{2q}\times\R$ endowed with the product
$$\langle(u,v),(u,v)\rangle=\langle u,u\rangle_q+v.v.$$We denote by $(e_1,\bar{e}_1,\ldots,e_q,\bar{e}_q)$ and $f$  the canonical basis of $\R^{2q}$ and $\R$, respectively.
We consider the skew-symmetric endomorphism $J$ of $V$ given by
 \begin{eqnarray*}
 Je_1&=&0,\;
 J\bar{e}_1=e_2,\;Jf=\bar{e}_2\\
 Je_2&=& \bar{e}_3-f,\;
 J\bar{e}_2=-e_1+\be e_3,\\
 Je_3&=&-\bar{e}_2+a_3e_3,\;
 J\bar{e}_3=-\be e_2-a_3\bar{e_3},\\
 Je_j&=&a_je_j,\;
 J\bar{e}_j=-a_j\bar{e}_j, j=4,\ldots,q.
 \end{eqnarray*}where  $0<a_{3}\leq\ldots\leq a_q$ and $\di\sum_{i=3}^qa_i^2=2\be.$
 One can check easily that $\ker J=\R e_1$ and ${\mathrm tr}(J^2)=0$. Thus  the bracket $[u,v]=\langle Ju,v\rangle e_1$ induces on $V$ a Heisenberg's Lie algebra structure for which $\prs$ is Ricci flat.
Put
 \begin{eqnarray*}
 E&=&e_1,\; F_1=f,\; \bar{F}_1={e}_2,\;F_2=\bar{e}_1,\;\bar{F}_2=\bar{e}_2,\\
 F_3&=&\frac1{\sqrt{a_3}}(f+e_3),\;\bar{F}_3=\frac1{\sqrt{a_3}}(\be\bar{e}_1+\bar{e}_3),\\
 F_i&=&\frac1{\sqrt{a_i}}e_i,\;\bar{F}_i=\frac1{\sqrt{a_i}}\bar{e}_i,\; i=4,\ldots,q.
 \end{eqnarray*}
 $(E,\bar{F}_1,\ldots,F_q,\bar{F}_q)$
is a basis of $V$ satisfying
$$\label{bracketexam}[F_j,\bar{F}_j]=E,\;\; j=1,\ldots,q,$$ and the all the others brackets are null or given by anti-symmetry. Thanks to this basis we identify $(V,[\;,\;],\prs)$ with the Lie algebra
$${\h}_{2q+1}=\left\{\left(\begin{array}{ccc}0&X&z\\0&0&{}^t\bar{X}\\0&0&0\end{array}\right),z\in\R,X,\bar{X}\in
\R^{q}\right\},$$endowed with the pseudo-Euclidean product of signature $(q,q+1)$ whose expression in the canonical basis $(e,m_j,\bar{m}_j)$  $j=1,\ldots,q$ of
${\h}_{2q+1}$ is given by
\begin{eqnarray*}
\langle e,m_2 \rangle&=&\langle m_1,m_1 \rangle=\langle \bar{m}_1,\bar{m}_2 \rangle=1,\\
\langle e, \bar{m}_3\rangle&=&\frac{\be}{\sqrt{a_3}},\; \langle m_3,m_3 \rangle=a_3^{-1},\;\langle m_1,m_3 \rangle=\frac{1}{\sqrt{a_3}},\\
\langle m_i,\bar{m}_i \rangle&=&a_i^{-1},\; i=3,\ldots,q.\end{eqnarray*}
The simply connected Lie group associated to ${\h}_{2q+1}$ is given by
$${H}_{2q+1}=\left\{\left(\begin{array}{ccc}1&X&z\\0&I&{}^t\bar{X}\\0&0&1\end{array}\right),z\in\R,X,\bar{X}\in\R^{q}\right\},$$and if
$(z,x_i,\bar{x}_i)$ are its canonical coordinates then the left invariant pseudo-Riemannian metric $\prs^l$ on ${H}_{2q+1}$ associated to $\prs$ is given by
\begin{eqnarray*}
 \prs^l&=&2\left(d\bar{x}_2+\frac{\be}{\sqrt{a_3}}d\bar{x_3}\right)\left(dz-\sum_{i=1}^{q}x_id\bar{x}_i\right)
 +2d\bar{x}_1d\bar{x}_{2}\\
 &&+\left(dx_1+\frac{1}{\sqrt{a_3}}dx_3\right)^2+
 2\sum_{i=3}^{q}a_i^{-1}dx_id\bar{x_i},\end{eqnarray*}where  $0<a_{3}\leq\ldots\leq a_q$ and $\di\sum_{i=3}^qa_i^2=2\be.$

 \end{enumerate}

\end{exem}

\section{Ricci flat Lorentzian 2-step nilpotent Lie groups}\label{section4}
It is well-known that a Ricci flat Euclidean Lie algebra must be flat (see \cite{besse}). In  \cite{milnor}, Milnor showed that an Euclidean Lie algebra is flat if and only if it is a semi-direct product of an abelian algebra $\mathfrak{b}$ with an abelian ideal $\mathfrak{u}$ and, for any $u\in\mathfrak{b}$, $ad_u$ is skew-symmetric. Thus an Euclidean 2-step nilpotent Lie algebra cannot be flat. In \cite{guediri2}, Guediri showed that a Lorentzian 2-step nilpotent Lie algebra is flat is and only if it is a trivial central extension of the 3-dimensional Heisenberg Lie algebra.
In this section, we determine all Ricci flat Lorentzian 2-step nilpotent Lie algebras.
\begin{theo}\label{lorentzflat}
Let $\mathfrak{N}$ be an irreducible Lorentzian 2-step nilpotent Lie algebra. Then $\mathfrak{N}$ is Ricci flat if and only if the center of $\mathfrak{N}$ is degenerate and  $\mathfrak{N}$ is isomorphic to  the Lorentzian vector space $\R^{(1,1)}\times\R^{p}\times\R^{2r}\times\R^q$ and, if $$(e,\bar e),\;(f_1,\ldots,f_p),\;(g_1,\ldots,g_{2r}),\;(h_1,\ldots,h_q)$$ are the canonical basis, respectively, of $\R^{(1,1)}$, $\R^{p}$, $\R^{2r}$ and $\R^q$, then
    the  Lie brackets are
    \begin{eqnarray*}
    \;[\bar e,g_i]&=&a_ie+\sum_{l=1}^px_i^lf_l,\; i=1,\ldots,2r,\\
    \;[\bar e,h_i]&=&b_ie+\sum_{l=1}^py_i^lf_l,\; i=1,\ldots,q,\\
    \;[g_{2i-1},g_{2i}]&=&\la_i e,\; i=1,\ldots,r,\end{eqnarray*} the others vanish or are obtained by symmetry, and the structure coefficients satisfy the following conditions:
    \begin{enumerate}\item $span\{(b_1,\ldots,b_q),(y_1^1,\ldots,y_q^1),\ldots,(y_1^p,\ldots,y_q^p)\}=\R^q,$
    \item $0<\la_1\leq\ldots\leq\la_r$ and
    $$\sum_{i,l}(x_i^l)^2+\sum_{i,l}(y_i^l)^2=\sum_{i=1}^r\la_i^2.$$

\end{enumerate}\end{theo}

{\it Proof.} Suppose that  $\mathfrak{N}$ is a Ricci flat irreducible Lorentzian 2-step nilpotent Lie algebra. We distinguish three cases:
\begin{enumerate}\item The center $\mathfrak{Z}$ of $\mathfrak{N}$ is nondegenerate and the restriction of $\prs$ to it is Lorentzian. In this case, according to (\ref{ricci}) and to the fact that $\mathfrak{Z}\subset \ker{\cal J}^-$ and $\mathfrak{Z}^\perp\subset \ker{\cal J}^+$, the vanishing of $\mathfrak{r}$ is equivalent to ${\cal J}^+={\cal J}^-=0.$\\
Choose a basis $\B=(e,\bar e,f_1,\ldots,f_{p},g_1,\ldots,g_{q})$ of $\mathfrak{N}$ such that, $(e,\bar e,f_1,\ldots,f_{p})$ is a Lorentzian basis of $\mathfrak{Z}$ and  $(g_1,\ldots,g_{q})$ is an orthonormal basis of $\mathfrak{Z}^\perp$, and denote by $(K, \bar K,J_1,\ldots,J_{p})$ the structure endomorphisms associated to $(e,\bar e,f_1,\ldots,f_{p})$.

Since $\mathfrak{Z}=\di \ker K\cap  \ker  \bar K\bigcap_{i=1}^{p}\ker J_i$, the representations of $K,\bar K,J_1,\ldots,J_{p}$ in $\B$ must have, respectively, the forms
$(0,M_1,0,0)$, $(0,M_2,0,0)$, $(0,B_i,0,0)$, $i=1,\ldots,p$. On the other hand, by using (\ref{invariant}), we get
$$\langle{\cal J}^+e,e\rangle=-\frac14{\mathrm tr}(\bar K^2),\;
\langle{\cal J}^+\bar e,\bar e\rangle=-\frac14{\mathrm tr}(K^2)\;\mbox{and}\;
\langle{\cal J}^+f_i,f_i\rangle=-\frac14{\mathrm tr}(J_i^2),$$for $i=1,\ldots,p$.
Thus ${\cal J}^+=0$ implies $${\mathrm tr}(K^2)={\mathrm tr}(\bar K^2)={\mathrm tr}(J_i^2)=0, \;i=1,\ldots,p,$$ and by using (\ref{trace}), we get $${\mathrm tr}(M_1^2)={\mathrm tr}(M_2^2)={\mathrm tr}(B_i^2)=0, \;i=1,\ldots,p.$$ Thus $K=\bar K=J_i=0$, for $i=1,\ldots,p$. So this case is impossible.

\item The center is nondegenerate and the restriction of $\prs$ to it is Euclidean. In this case, according to (\ref{ricci}) and to the fact that $\mathfrak{Z}\subset \ker{\cal J}^-$ and $\mathfrak{Z}^\perp\subset \ker{\cal J}^+$, the vanishing of $\mathfrak{r}$ is equivalent to ${\cal J}^+={\cal J}^-=0.$\\
     Choose a basis $\B=(e,\bar e,f_1,\ldots,f_{p},g_1,\ldots,g_{q})$ of $\mathfrak{N}$ such that, $(f_1,\ldots,f_{p})$ is an orthonormal  basis of $\mathfrak{Z}$ and  $(e,\bar e,g_1,\ldots,g_{q})$ is a Lorentzian basis of $\mathfrak{Z}^\perp$. Denote by $(J_1,\ldots,J_{p})$ the structure endomorphisms associated to $(f_1,\ldots,f_{p})$. By using (\ref{invariant}), we get
\begin{eqnarray*}
{\cal J}^+&=&-\frac14\sum_{i,j}\langle f_i,.\rangle{\mathrm tr}(J_i\circ J_j)f_j\quad
\mbox{and}\quad{\cal J}^-=\frac12\sum_{i=1}^{p}J_i^2.\end{eqnarray*}
For $i=1,\ldots,p$, let $(A_i,B_i,X^i,Y^i)$ denotes the representation of $J_i$ in $\B$ where $A_i=\left(\begin{array}{cc}a_i&0\\0&-a_i\end{array}\right)$.   Now, by using (\ref{product}) and (\ref{formule2}), we get that $\sum_{i=1}^{p}J_i^2=0$ implies that
$$\sum_{i=1}^p\left(\begin{array}{cc}
a_i^2-X^i.Y^i&-X^i.X^i\\-Y^i.Y^i&a_i^2-X^i.Y^i\end{array}\right)=0,$$
 and then $X^i=Y^i=a_i=0$ for $i=1,\ldots,p$. On the other hand ${\cal J}^+=0$ implies $
{\mathrm tr}(J_i^2)={\mathrm tr}(B_i^2)=0$ and hence $B_i=0$. Finally, $J_1=\ldots=J_p=0$ and we conclude that this case is impossible.

\item The center is degenerate. Choose a Lorentzian basis $$\B=(e,\bar e,f_1,\ldots,f_{p},g_1,\ldots,g_{s})$$ such that $(e,f_1,\ldots,f_{p})$ is a   basis of $\mathfrak{Z}$ and  $(e, g_1,\ldots,g_{s})$ is a  basis of $\mathfrak{Z}^\perp$ and denote by $(K,J_1,\ldots,J_{p})$ the structure endomorphisms associated to $(e,f_1,\ldots,f_{p})$.  By using (\ref{invariant}), we get
    \begin{eqnarray*}
    {\cal J}^-&=&\frac12\sum_{i=1}^{p}J_i^2,\\
    {\cal J}^+e&=&{\cal J}^+g_i=0,\; i=1,\ldots,s,\\
    {\cal J}^+\bar e&=&-\frac14{\mathrm tr}(K^2)e-\frac14\sum_{i=1}^{p}
    {\mathrm tr}(K\circ J_i)f_i,\\
    {\cal J}^+f_j&=&-\frac14\sum_{i=1}^{p}{\mathrm tr}(J_j\circ J_i)f_i,\; j=1,\ldots,p.\end{eqnarray*}So we can deduce from (\ref{ricci}) that the vanishing of $\mathfrak{r}$ is equivalent to
    \begin{equation}\label{temp1}\left\{\begin{array}{ccl}
    {\mathrm tr}(J_j\circ J_i)&=&{\mathrm tr}(K\circ J_i)=0,\; i,j=1,\ldots,p,\\{\cal J}^-&=&\frac12\langle e,.\rangle{\mathrm tr}(K^2)e.\end{array}\right.\end{equation}
    Since $\di\mathfrak{Z}= \ker K\cap\bigcap_{i=1}^{p}\ker J_i$, the representation of
    $(K,J_1,\ldots,J_{p})$ in $\B$ are given by $((0,M,V,0),(0,B_1,X^1,0),\ldots,
    (0,B_{p},X^{p},0))$ and, for $i=1,\ldots,p$,
    \begin{eqnarray*}
    V&=&(0,\ldots,0,a_1,\ldots,a_{s}),\; X^i=(0,\ldots,0,x_1^i,\ldots,x_{s}^i),\\
    M&=&\left(\begin{array}{cc}0&0\\0&M_0\end{array}\right),\;
    B_i=\left(\begin{array}{cc}0&0\\0&B_{0,i}\end{array}\right),\end{eqnarray*}where $M_0$ and $B_{0,i}$ are skew-symmetric $(s,s)$ matrix.

    We can suppose that $M_0$ has the form (\ref{diagonal}).\\
     Now, we have ${\mathrm tr}J_i^2={\mathrm tr}B_{0,i}^2=0$ and hence $B_1=\ldots=B_{p}=0.$
On the other hand, by using (\ref{product})-(\ref{formule2}) we can see easily that  ${\cal J}^-=\frac12\langle e,.\rangle{\mathrm tr}(K^2)e$ is equivalent to
\[\di\sum_{i=1}^{p}X^i.X^i=-\frac12{\mathrm tr}(K^2)=\sum_{j=1}^r\la_j^2.\]
To summarize, we have that (\ref{temp1}) is equivalent to
\begin{equation}
\label{temp2}B_1=\ldots=B_{p}=0\quad\mbox{and}\quad\sum_{i=1}^{p}X^i.X^i
=\sum_{j=1}^r\la_j^2.\end{equation}

Let us find the necessary and the sufficient conditions to have $\di\mathfrak{Z}= \ker K\cap\bigcap_{i=1}^{p}\ker J_i$.\\
A vector of coordinates $(a,b,z_1,\ldots,z_{p},t_1,\ldots,t_{q})$ is an element of $\ker K\cap\bigcap_{i=1}^{p}\ker J_i$ if and only if
$$\left\{\begin{array}{l}a_1t_1+\ldots+a_{s}t_{s}=0,\\
x_1^it_1+\ldots+x_{s}^it_{s}=0,\; i=1,\ldots,p,\\
a_{2i-1}b-\la_it_{2i}=a_{2i}b+\la_it_{2i-1}=0,i=1,\ldots,r\\
a_jb=bx_k^i=0, j=2r+1,\ldots,s, i=1,\ldots,p,k=1,\ldots,s.\end{array}\right.$$
We have two cases:
\begin{enumerate}\item $X_1=\ldots=X_{p}=0$. In this case the condition (\ref{temp2}) implies $\la_1=\ldots=\la_r=0$ and the system above reduces to
$$\left\{\begin{array}{l}a_1t_1+\ldots+a_{s}t_{s}=0,\\
a_jb=0, j=1,\ldots,s.\end{array}\right.$$
The condition $\ker K=\mathfrak{Z}$ is equivalent to $s=1$ and $a_1\not=0$. Thus in the basis $(e,\bar e,f_1,\ldots,f_{p},g_1)$ the non vanishing brackets are given by
$$[\bar e,g_1]=ae.$$Hence $\mathfrak{N}$ is a trivial central extension of the 3-dimensional Heisenberg Lie algebra and the metric is flat.
\item There exists $i$ such that $X_i\not=0$. In this case the system is equivalent to
$$\left\{\begin{array}{l}a_{2r+1}t_{2r+1}+\ldots+a_{s}t_{s}=0,\\
x_{2r+1}^it_{2r+1}+\ldots+x_{s}^it_{s}=0,\; i=1,\ldots,p,\\
b=t_1=\ldots=t_{2r}=0,\\
\end{array}\right.$$
and the condition $\ker K\cap\bigcap_{i=1}^{p}\ker J_i=\mathfrak{Z}$ if and only if the subspace
$$span\left\{(a_{2r+1},\ldots,a_{s}),(x_{2r+1}^1,\ldots,x_{s}^1),\ldots,
(x_{2r+1}^{p},\ldots,x_{s}^{p})\right\}$$ of $\R^{s-2r}$ is of dimension $q=s-2r$.
By putting $(y_1^l,\ldots,y_q^l)=(x_{2r+1}^l,\ldots,x_{s}^l)$, $l=1,\ldots,p$ and $(b_1,\ldots,b_q)=(a_{2r+1},\ldots,a_s)$, we get the desired result. $\;\square$

\end{enumerate}
\end{enumerate}

 From Theorem \ref{lorentzflat} and by a direct computation we get  the following result which give all Ricci flat left invariant Lorentzian metrics on 2-step nilpotent Lie groups.
\begin{theo}\label{main}Let $(N,\prs)$ be a 2-step nilpotent Lie group endowed with a left invariant Lorentzian metric. Then $\prs$ is Ricci flat if and only if:
 \begin{enumerate}\item there exists $M_1=(x_i^j)$ and $M_2=(y_i^j)$, respectively, a real $(2r,p)$ matrix and $(q,p)$ matrix, $0<\lambda_1\leq\ldots\leq\lambda_r$, $A\in\R^{2r}$ and $B\in\R^q$ such that
 $${\mathrm tr}(M_1{}^tM_1)+{\mathrm tr}(M_2{}^tM_2)=\sum_{i=1}^r\la_i^2,$$
 and $span\{B,Y_1,\ldots,Y_p\}=\R^q$, ($Y_1,\ldots,Y_p$ are column's vectors of $M_2$),
\item $(N,\prs)$ is isomorphic to $(\R^2\times\R^{p}\times\R^{2r}\times\R^q,\prs_0)$ where:
 \begin{enumerate} \item the Lie group structure is given by
 \begin{eqnarray*}
&&(t_1,\bar{t}_1,U_1,V_1,W_1).(t_2,\bar{t}_2,U_2,V_2,W_2)=(T,U,V,W),\\
T&=&(t_1+t_2+
\frac12\sum_{i=1}^r\la_i\left(v_{2i-1}^1v_{2i}^2-v_{2i-1}^2v_{2i}^1\right)
+\frac12\sum_{i=1}^{2r}a_i\left(\bar{t}_1v_{i}^2-\bar{t}_2v_{i}^1\right)\\&&
+\frac12\sum_{i=1}^qb_i\left(\bar{t}_1w_{i}^2-\bar{t}_2w_{i}^1\right),\bar{t}_1+\bar{t_2}),\\
U&=&\left(u_1^1+u_1^2+\frac12\sum_{i=1}^{2r}x_i^1\left(\bar{t}_1v_{i}^2-\bar{t}_2v_{i}^1\right)
+\frac12\sum_{i=1}^{q}y_i^1\left(\bar{t}_1w_{i}^2-\bar{t}_2w_{i}^1\right),\ldots,\right.\\&&\left.
u_p^1+u_p^2+\frac12\sum_{i=1}^{2r}x_i^p\left(\bar{t}_1v_{i}^2-\bar{t}_2v_{i}^1\right)
+\frac12\sum_{i=1}^{q}y_i^p\left(\bar{t}_1w_{i}^2-\bar{t}_2w_{i}^1\right)\right)\\
V&=&(v_1^1+v_1^2,\ldots,v_{2r}^1+v_{2r}^2),\;
W=(w_1^1+w_1^2,\ldots,w_{q}^1+w_{q}^2).
\end{eqnarray*}
\item in the canonical coordinates $(t,\bar{t},u_1,\ldots,u_p,v_1,\ldots,v_{2r},w_1,\ldots,w_q)$ the metric is given by (the omitted products are null):
\begin{eqnarray*}
\langle \partial_{t},\partial_{\bar{t}} \rangle_0&=&\langle \partial_{u_i},\partial_{u_i} \rangle_0=1,\;
\langle \partial_{\bar{t}},\partial_{u_i} \rangle_0=\frac12\left(X_i.V+Y_i.W\right),\;i=1,\ldots,p,\\
\langle \partial_{\bar{t}},\partial_{\bar{t}} \rangle_0&=&A.V+B.W+\frac14\sum_{l=1}^p\left(X_l.V+Y_l.W\right)^2,\\
\langle \partial_{u_l}, \partial_{v_i}\rangle_0&=&-\frac{\bar{t}}2x_{i}^l,\;
\langle \partial_{u_l}, \partial_{w_j} \rangle_0=-\frac{\bar{t}}2y_{j}^l,\;l=1,\ldots,p,\;i=1,\ldots,2r,\;j=1,\ldots,q,\\
\langle \partial_{\bar{t}},\partial_{v_i} \rangle_0&=&-\frac12((-1)^i\la_{[\frac{i+1}2]}v_{i-(-1)^i}+a_{i}\bar{t})-
\frac{\bar{t}}4X^i.\left(\sum_{j=1}^{2r}v_jX^j+\sum_{j=1}^{q}w_jY^j\right),\\&&\;i=1,\ldots,2r,\\
\langle \partial_{\bar{t}},\partial_{w_i} \rangle_0&=&-\frac{\bar{t}}2b_i-
\frac{\bar{t}}4Y^i.\left(\sum_{j=1}^{2r}v_jX^j+\sum_{j=1}^{q}w_jY^j\right),\;i=1,\ldots,q,\\
\langle \partial_{v_i} , \partial_{v_j}\rangle_0&=&\de_{ij}+\frac{\bar{t}^2}4X^i.X^j,\;i,j=1,\ldots,2r,\\
\langle \partial_{v_i} , \partial_{w_j} \rangle_0&=&\frac{\bar{t}^2}4X^i.Y^j,\;i=1,\ldots,2r,j=1,\ldots,q,\\
\langle \partial_{w_i} , \partial_{w_j} \rangle_0&=&\de_{ij}+\frac{\bar{t}^2}4Y^i.Y^j,\;i,j=1,\ldots,q,\\
\end{eqnarray*}where $(X_1,\ldots,X_{p})$ (resp. $(X^1,\ldots,X^{2r})$) are the columns's vectors (resp. the row's vectors) of $M_1$ and $(Y_1,\ldots,Y_{p})$ (resp. $(Y^1,\ldots,Y^{q})$) are the columns's vectors (resp. the row's vectors) of $M_2$.

\end{enumerate}\end{enumerate}
\end{theo}
\begin{rem}\begin{enumerate}\item For $p=r=0$ and $q=1$ we get the lowest dimensional 2-step nilpotent Lie group with Ricci flat Lorentzian metric. It is isomorphic to $\R^3$ with the product
$$(t_1,\bar{t}_1,w_1).(t_2,\bar{t}_2,w_2)=(t_1+t_2+\frac12(\bar{t}_1w_2-\bar{t}_2w_1),\bar{t}_1+\bar{t}_2,w_1+w_2)$$
and the metric is given by
$$\prs_0=dw^2+d\bar{t}\left(wd\bar{t}-\bar{t}dw+2dt\right).$$
We recover the 3-dimensional Heisenberg group with its canonical flat Lorentzian metric.
\item There is no irreducible 4-dimensional 2-step nilpotent Lie group with Ricci flat Lorentzian metric and for $p=r=1$ and $q=0$ we get the second low dimensional 2-step nilpotent Lie group with Ricci flat Lorentzian metric. It is isomorphic to $\R^5$ with the product
$$(t_1,\bar{t}_1,u_1,v_1^1,v_2^1).(t_2,\bar{t}_2,u_2,v_1^2,v_2^2)=(T,\bar{t}_1+\bar{t}_2,U,v_1^1+v_1^2,v_2^1+v_2^2)$$
where
\begin{eqnarray*}
T&=&t_1+t_2+
\frac12\sqrt{\al^2+\be^2}\left(v_{1}^1v_{2}^2-v_{1}^2v_{2}^1\right)
\\&&+\frac12\left(\bar{t}_1(a_1v_{1}^2+
a_2v_{2}^2)-\bar{t}_2(a_1v_{1}^1+
a_2v_{2}^1)\right)\\
U&=&u_1+u_2+\frac12\left(\bar{t}_1(\al v_{1}^2+
\be v_{2}^2)-\bar{t}_2(\al v_{1}^1+
\be v_{2}^1)\right).
\end{eqnarray*} ($\al,\be,a_1,a_2$ are parameters satisfying $(\al,\be)\not=(0,0)$).
The metric is given by (the omitted products are null):
\begin{eqnarray*}
\langle \partial_{t},\partial_{\bar{t}} \rangle_0&=&\langle \partial_{u},\partial_{u} \rangle_0=1,\;
\langle \partial_{\bar{t}},\partial_{u} \rangle_0=\frac12(\al v_1+\be v_2),\\
\langle \partial_{\bar{t}},\partial_{\bar{t}} \rangle_0&=&a_1 v_1+a_2 v_2+\frac14(\al v_1+\be v_2)^2,\\
\langle \partial_{u}, \partial_{v_1}\rangle_0&=&-\frac{\bar{t}}2\al,\;
\langle \partial_{u}, \partial_{v_2} \rangle_0=-\frac{\bar{t}}2\be,\\
\langle \partial_{\bar{t}},\partial_{v_1} \rangle_0&=&-\frac12(-\sqrt{\al^2+\be^2}v_{2}+a_{1}\bar{t})-
\frac{\bar{t}}4\al(v_1\al+v_2\be),\\
\langle \partial_{\bar{t}},\partial_{v_2} \rangle_0&=&-\frac12(\sqrt{\al^2+\be^2}v_{1}+a_{2}\bar{t})-
\frac{\bar{t}}4\be(v_1\al+v_2\be),\\
\langle \partial_{v_1} , \partial_{v_1}\rangle_0&=&1+\frac{\bar{t}^2}4\al^2,\;
\langle \partial_{v_1} , \partial_{v_2} \rangle_0=\frac{\bar{t}^2}4\al\be,\;
\langle \partial_{v_2} , \partial_{v_2} \rangle_0=1+\frac{\bar{t}^2}4\be^2.
\end{eqnarray*}

\end{enumerate}

\end{rem}

\paragraph{Acknowledgement:}

 A part of this work was done at The Abdus Salam Centre of Theoretical Physics, the author would like to thank the Mathematic section for hospitality.\\

{\it Mohamed BOUCETTA,\\ Facult\'e des sciences et techniques Gueliz\\BP 549 Marrakech Maroc\\
mboucetta2@yahoo.fr}

\end{document}